\documentclass[hidelinks]{article}
\usepackage{graphicx} 
\usepackage{hyperref}
\usepackage{amsthm}
\usepackage{geometry}
\usepackage{enumerate}
\usepackage{amssymb}
\usepackage{amsmath}
\usepackage{xcolor}
\definecolor{purple}{rgb}{0.5,0,0.5}
\usepackage{lmodern}
\usepackage[english]{babel}

\newtheorem{theorem}{Theorem}[section]
\newtheorem{lemma}[theorem]{Lemma}

\newtheorem{claim}[theorem]{Claim}

\theoremstyle{definition}
\newtheorem{definition}[theorem]{Definition}

\newcommand{\claimproofstart}[1][Proof]{\begin{proof}[#1]
\renewcommand{\qedsymbol}{$\boxdot$}}
\newcommand\claimproofend{
\end{proof}
\renewcommand{\qedsymbol}{$\square$}}

\newcommand{\cF}{\mathcal{F}}

\newcommand{\cH}{\mathcal{H}}

\newcommand{\rind}{r_{\mathrm{ind}}}
\newcommand{\ramseyind}{\xrightarrow[]{\text{ind}}}

\title{Erd\H{o}s--Szekeres bound for induced Ramsey numbers}

\author{Rajko Nenadov \footnote{School of Mathematics and Statistics, University of Canterbury, New Zealand. Research supported by the New Zealand Marsden Fund. Email: rajko.nenadov@canterbury.ac.nz}}

\begin{document}

\maketitle

\begin{abstract}
 The Erd\H{o}s--Szekeres bound from 1935 gives $\log r(F,H)=O(s\log(2t/s))$ for all graphs $F$ and $H$ on $s$ and $t$ vertices, respectively, where $2\le s\le t$. We prove that this bound extends to induced Ramsey numbers. The dependence on $s$ and $t$ is optimal in general, up to an absolute constant factor.
\end{abstract}

\section{Introduction}

Ramsey's theorem~\cite{Ram30} asserts that for every two graphs $F$ and $H$, any red/blue colouring of the edges of a sufficiently large complete graph $K_n$ contains a blue copy of $F$ or a red copy of $H$. The smallest such number $n$, denoted by $r(F,H)$, is called the \emph{Ramsey number}. In 1935, Erd\H{o}s and Szekeres \cite{erdos35szekeres} proved
$$
    r(F,H) \le \binom{s+t-2}{s-1},
$$
where $s=v(F)$ and $t=v(H)$. For $2\le s\le t$, this gives
\begin{equation}\label{eq:classical-ramsey}
    \log r(F,H)=O\!\left(s\log\frac{2t}{s}\right).
\end{equation}
A standard probabilistic argument, following Erd\H{o}s~\cite{Erdos1947}, gives
$$
     \log r(K_s,K_t)=\Omega\!\left(s\log\frac{2t}{s}\right),
$$
for $2\le s\le t$; thus the dependence on $s$ and $t$ in \eqref{eq:classical-ramsey} is optimal up to an absolute constant factor. Obtaining sharper estimates, particularly for complete graphs, has been a central problem in combinatorics (see, e.g., \cite{BBCGHMST2026,CGMS2026,Conlon2009,HeftyHornKingPfender2025,MaShenXie2026,MattheusVerstraete2024,Sah2023}).

Requiring the monochromatic copies to be induced introduces a further constraint: nonedges of the target graph must remain nonedges of the host. Formally, we say that a graph $G$ is \emph{induced Ramsey} for the pair of graphs $F$ and $H$, and write
$$
    G \ramseyind (F,H),
$$
if in every red/blue colouring of the edges of $G$, there exists an induced copy of $F$ in $G$ with all its edges blue, or an induced copy of $H$ in $G$ with all its edges red. It was proved independently by Deuber \cite{Deuber1975}, Erd\H{o}s, Hajnal, and P\'osa \cite{ErdosHajnalPosa1975}, and R\"odl \cite{Rodl1973} that such a graph $G$ exists for every pair of graphs $F$ and $H$. The \emph{induced Ramsey number} $\rind(F,H)$ is the minimum number of vertices in a graph that is induced Ramsey for $(F,H)$.

Building on work by Kohayakawa, Pr\"omel, and R\"odl \cite{KohayakawaPromelRodl1998},  Conlon, Fox, and Sudakov \cite{ConlonFoxSudakov2012} proved that $\log \rind(F, F) = O(s \log s)$, where $s = v(F)$. In a recent breakthrough, Arag{\~a}o, Campos, Dahia, Filipe, and Marciano \cite{AragaoCamposDahiaFilipeMarciano2025} proved that
$$
    \log \rind(F,F) \le C s,
$$
where $C>0$ is an absolute constant. This resolved a long-standing conjecture of Erd\H{o}s \cite{Erdos1984Problems}, matching the Erd\H{o}s--Szekeres bound up to an absolute constant factor in the exponent.

In the general, off-diagonal setting, Kohayakawa, Pr\"omel, and R\"odl~\cite{KohayakawaPromelRodl1998} proved that $\log\rind(F,H)=O(s(\log t)^2)$, where $s=v(F)$ and $t=v(H)$. They conjectured that, for each fixed graph $F$, the induced Ramsey number $\rind(F,H)$ is polynomial in $t$. This conjecture was resolved independently by Zhong, Kashima, and Mao \cite{ZhongKashimaMao2026} and by Andrades, Campos, Filipe, and Kats \cite{AndradesCamposFilipeKats2026}, who proved
$$
    \log \rind(F,H) = O(s^2 \log t).
$$
Applying the diagonal bound to the disjoint union $F\sqcup H$ gives $\log\rind(F,H)=O(t)$ when $s\le t$. Together, these results yield
$$
    \log\rind(F,H)=O\!\left(\min\{t,s^2\log t,s(\log t)^2\}\right)
$$
for $2\le s\le t$.

We prove that the Erd\H{o}s--Szekeres upper bound~\eqref{eq:classical-ramsey} extends to induced Ramsey numbers throughout the full range $2\le s\le t$.

\begin{theorem} \label{thm:main}
    There exists an absolute constant $C > 1$ such that for all graphs $F$ and $H$ on $s$ and $t$ vertices, respectively, where $1\le s\le t$, we have
    $$
        \log \rind(F, H) \le Cs \log \frac{2 t}{s}.
    $$
\end{theorem}

For complete graphs, induced and ordinary Ramsey numbers coincide: $\rind(K_s,K_t)=r(K_s,K_t)$. Together with the lower bound above, this shows that the dependence on $s$ and $t$ in Theorem~\ref{thm:main} is optimal in general, up to an absolute constant factor. Taking $F=H$ recovers the two-colour bound of Arag{\~a}o et~al.~\cite{AragaoCamposDahiaFilipeMarciano2025}. The theorem also gives $\rind(F,H)\le t^{O(s)}$, answering the question of Andrades, Campos, Filipe, and Kats~\cite[Question~6.1]{AndradesCamposFilipeKats2026}. The dependence on the ratio $t/s$ also improves the bound in the intermediate range. For example, when $t=s^2$, the preceding results give $\log\rind(F,H)=O(s(\log s)^2)$, whereas Theorem~\ref{thm:main} gives $O(s\log s)$.

Stronger bounds are known for restricted classes of $t$-vertex graphs $F$ and $H$: polynomial bounds in $t$ for graphs with fixed maximum degree~\cite{ConlonFoxZhao2014,LuczakRodl1996}, and linear bounds in $t$ for cycles~\cite{HaxellKohayakawaLuczak1995} and, more generally, for graphs with both maximum degree and treewidth bounded by fixed constants~\cite{HunterSudakov2025}.

\paragraph{Proof overview.} For simplicity, suppose that $F$ has an even number of vertices, and let $V(F)=A\cup B$ be a partition with $|A| = |B|$. Suppose that we have constructed an $n$-vertex graph $G$ with a proper colouring into $\ell$ equally sized classes, such that
$$
    G[S] \ramseyind (F[A], H) \quad \text{ and } \quad G[S] \ramseyind (F[B], H)
$$
for every $S\subseteq V(G)$ of size at least $(1-1/(2t))n$. Write $V(H)=[t]$.

Using Lemma~\ref{lemma:robust}, we amplify $G$ to a graph $\Gamma$ with $m$ vertices and maximum degree at most $m^{0.98}$, such that the same Ramsey property holds on every subset of size at least $\beta m$, for a suitably small $\beta>0$. We now construct $G'$ on disjoint vertex sets $V_1,\ldots,V_t$, each of size $m$, as follows:
\begin{itemize}
    \item Put a distinct copy of $\Gamma$ inside each $V_i$, so that $G'[V_i] \cong \Gamma$.
    \item For $\{i,j\}\in E(H)$, include each pair between $V_i$ and $V_j$ as an edge independently with probability $1 - q$, where $q = c s / t$ for a suitably small constant $c > 0$.
    \item For $\{i,j\}\notin E(H)$, put no edges between $V_i$ and $V_j$.
\end{itemize}
With positive probability, every pair of sufficiently large subsets in parts corresponding to an edge of $H$ has density at least $1-2q$ in $G'$. Consider a red/blue colouring of the edges of $G'$ with no red induced copy of $H$. The elementary greedy embedding lemma, Lemma~\ref{lemma:density_embedding}, implies that there exist large sets $X_i\subseteq V_i$ and $X_j\subseteq V_j$, for some $\{i,j\}\in E(H)$, between which the red density is less than $1-3q$. Otherwise, Lemma~\ref{lemma:density_embedding} gives a red copy of $H$ meeting each part once; this copy is induced because parts corresponding to nonedges of $H$ have no edges between them. The density of $G'$ between $X_i$ and $X_j$ then gives blue density at least $q$. By the assumption on $\Gamma$, every subset of $X_i$ of size at least $\beta m$ contains an induced blue copy of $F[A]$, and the analogous statement holds for $F[B]$ in $X_j$.

The key probabilistic statement is that, with positive probability, the sampled graph has the following property simultaneously for every red/blue colouring with no red induced copy of $H$: some pair of blue induced copies of $F[A]$ and $F[B]$ is joined by exactly the prescribed pattern of blue edges and absent edges. Together, these copies form a blue induced copy of $F$. Lemma~\ref{lemma:join} establishes this, using an argument based on the Sauer--Shelah inequality which controls all possible blue crossing subgraphs simultaneously, avoiding a union bound over their edge sets (see Section \ref{sec:join}). The remaining union bounds are controlled by the maximum degree of $\Gamma$. 

The same argument applies after deleting at most a $1/(2t)$-fraction of the vertices of $G'$, as this implies that at least half the vertices of each $V_i$ remain, thus preserving the robustness required for the next step. Joining equal-sized pieces allows the target size to double at each step, and tracking the resulting recurrence gives the exponent $O(s\log(2t/s))$.

\paragraph{Comparison with previous work.} Our proof is closest in spirit to the diagonal argument of Arag{\~a}o et~al.~\cite{AragaoCamposDahiaFilipeMarciano2025}. They maintain well-distributed families of induced copies, quantified by a Janson-type second-moment condition, and use hypergraph containers to extend these copies by one vertex. We join copies on two equally sized pieces. In our induction, the required abundance is expressed by the existence of the appropriate induced blue copy in every sufficiently large vertex subset, assuming that there is no red induced copy of $H$. Our joining step combines dependent random choice with a probabilistic corollary of the Sauer--Shelah inequality to control all possible blue crossing subgraphs simultaneously.

Our construction also shares a geometric ingredient with the recent off-diagonal proofs~\cite{AndradesCamposFilipeKats2026,ZhongKashimaMao2026}. Namely, the robust graph supplied by Lemma \ref{lemma:robust} is obtained by overlapping random blow-ups of induced subgraphs of $G$, so that every sufficiently large vertex subset contains an almost-spanning induced subgraph of the original graph. Such constructions have precedents in the work of Brown and R\"odl~\cite{BrownRodl1991} and Kohayakawa, Pr\"omel, and R\"odl~\cite{KohayakawaPromelRodl1998}.

\paragraph{AI Declaration.} The author and ChatGPT Astra had previously discussed applications of the Sauer--Shelah inequality to the max-cut problem in random induced subgraphs. Building on those discussions, Astra suggested using the inequality to control the union bound in this paper, providing the central new idea of the proof. Astra also assisted with literature searches, the refinement and simplification of arguments, the choice of constants, and the exposition. The author takes full responsibility for the mathematical arguments, references, and final manuscript.

\section{Robustness amplification}

Each inductive step first amplifies the robustness of a previously constructed host and then joins copies of the resulting graph using bipartite random graphs. This section establishes the amplification step. The following definition captures one of the main properties of graphs that we construct. 

\begin{definition}[Robust Ramsey property]
Let $H$ be a graph, let $b \ge 1$ be an integer, and let $c \in (0, 1]$. We say that a graph $G$ is \emph{$(H, b, c)$-Ramsey} if 
$$
    G[S] \ramseyind (F, H)
$$
for every $S \subseteq V(G)$ of size $|S| \ge c v(G)$ and every graph $F$ with at most $b$ vertices. 
\end{definition}

The following lemma amplifies $(H, b, 1 - \alpha)$-Ramsey graphs to $(H, b, \beta)$-Ramsey graphs, for $0<\alpha<1$ and $0<\beta\le1$, with a controlled additional cost.  We say that a graph $G$ is \emph{$\ell$-equicolourable} if it admits a proper colouring with $\ell$ colours in which all colour classes have the same size. Such a colouring is called an \emph{$\ell$-equicolouring}.

\begin{lemma} \label{lemma:robust}
There exists an absolute constant $C>1$ such that the following holds. Let $H$ be a graph, let $\ell,b\ge1$ be integers, and let $0<\alpha<1$. Suppose that $G$ is an $\ell$-equicolourable $(H,b,1-\alpha)$-Ramsey graph on $n$ vertices. For every $0< \beta \le 1$ and every prime power $r\ge d$, where
$$
    d=\left\lceil C\,\frac{\ell\log(2\ell)}{\alpha\beta^2}\right\rceil,
$$
there exists an $r$-equicolourable $(H,b,\beta)$-Ramsey graph $\Gamma$ with $m$ vertices, satisfying
$$
    m = \frac{r^2}{\ell}\, n \qquad\text{and}\qquad \Delta(\Gamma) \le \frac{d}{r}\,m.
$$
\end{lemma}

The proof of Lemma~\ref{lemma:robust} follows ideas of Brown and R\"odl~\cite{BrownRodl1991}. They showed that, for each fixed $\beta\in(0,1)$ and every $n$-vertex graph $G$, there exists a graph $\Gamma$ with $O_\beta(n^2\log^2 n)$ vertices such that every $S\subseteq V(\Gamma)$ with $|S|\ge\beta v(\Gamma)$ contains an induced copy of $G$. The order of this host is too large for our induction. Moreover, quadratic dependence on $n$ is necessary in the worst case over $n$-vertex target graphs. Our variant exploits an $\ell$-equicolouring of $G$ and its robustness under vertex deletion: it is enough to recover an induced subgraph of $G$ on at least $(1-\alpha)n$ vertices, which retains the required Ramsey property. The maximum-degree bound on $\Gamma$ will also be important when taking union bounds later in the proof.

Recall that a hypergraph is \emph{linear} if any two distinct hyperedges meet in at most one vertex. The hypergraph provided by the next lemma serves as a ``blueprint'' for overlapping blow-ups of $G$, which gives a robust graph $\Gamma$. 

\begin{lemma} \label{lemma:geometry}
    Let $0<\beta\le1$, let $d\ge16/\beta$ be an integer, and let $r\ge d$ be a prime power. Then there exists a linear $r$-uniform hypergraph $\cH$ with $r^2$ vertices and maximum degree $d$, with the following properties:
    \begin{enumerate}        
        \item \label{prop:transversal} there exists a partition $V(\cH) = X_1 \cup \cdots \cup X_r$ into sets of size $r$, such that every hyperedge meets each $X_i$ exactly once, and
        \item for every $A \subseteq V(\cH)$ of size $|A| \ge \beta r^2 / 2$, at least $d r /2$ hyperedges $h \in \cH$ satisfy
        \begin{equation}\label{eq:robust-rich-lines}
             |A\cap h|\ge\frac{|A|}{2r}.
        \end{equation}
    \end{enumerate}
\end{lemma}
\begin{proof}
Choose an arbitrary set $D \subseteq \mathbb F_r$ of size  $|D|=d$. We define $\cH$ on the vertex set $\mathbb F_r^2$, and for each $a \in D$ and $c \in \mathbb F_r$  take the hyperedge
$$
    h_{a,c}=\{(x,ax+c) \; :\; x \in\mathbb F_r\}.
$$
This hypergraph is clearly linear and $d$-regular. The desired partition is given by the vertical lines $\{(a,x):x\in\mathbb F_r\}$, indexed by $a\in\mathbb F_r$. It remains to verify the last property. 

Fix $A\subseteq\mathbb F_r^2$ of size $|A| \ge \beta r^2/2$, and let $\mathcal L$ be the set of all affine lines, including the vertical ones. Each point lies on exactly $r+1$ affine lines, and each ordered pair of distinct points lies on exactly one affine line. Thus
$$
    \sum_{L\in\mathcal L}|A\cap L|=(r+1)|A|,
    \quad \text{ and } \quad
    \sum_{L\in\mathcal L}|A\cap L|(|A\cap L|-1) =|A|(|A|-1).
$$
In particular,
$$
    \sum_{L\in\mathcal L}|A\cap L|^2=|A|^2+r|A|.
$$
Since $|\mathcal L|=r(r+1)$, expanding the square gives
$$
    \begin{aligned}
    \sum_{L\in\mathcal L}
    \left(|A\cap L|-\frac{|A|}{r}\right)^2
    &= |A|^2+r|A| - \frac{2(r+1)|A|^2}{r} + \frac{(r+1)|A|^2}{r}\\
    &= r|A|-\frac{|A|^2}{r} \le r|A|.
    \end{aligned}
$$
Every line failing~\eqref{eq:robust-rich-lines} contributes more than $|A|^2/(4r^2)$ to this sum. Thus at most
$$
    \frac{4r^3}{|A|}\le\frac{8r}{\beta}\le\frac{dr}{2}
$$
of the $dr$ chosen lines fail~\eqref{eq:robust-rich-lines}.
\end{proof}

We replace each vertex of $\cH$ by an independent set and place a random blow-up of an induced subgraph of $G$ on the sets indexed by each hyperedge. The incidence property and independence of the choices on different hyperedges ensure, with positive probability, that every large vertex set contains an induced copy of $G[T]$ for some $|T|\ge(1-\alpha)n$, and hence inherits the required Ramsey property.

\begin{proof}[Proof of Lemma \ref{lemma:robust}]
Let $V_1,\ldots,V_\ell$ be a partition of $V(G)$ corresponding to an $\ell$-equicolouring, so $|V_i|=n/\ell$. For sufficiently large $C$, we have $d\ge16/\beta$. Let $\cH$ be a hypergraph given by Lemma~\ref{lemma:geometry} with parameters $r,d,\beta$. The assumption $r \ge d$ gives the following inequality, which we shall use a couple of times: For every $A \subseteq V(\cH)$ of size $|A| \ge \beta r^2 / 2$, we have
\begin{equation} \label{eq:alpha_beta_l}
    \frac{\alpha \beta |A|}{r} \ge \frac{C}{2} \ell \log(2\ell).
\end{equation}

For each $x\in V(\cH)$, introduce a new set $C_x$ of size $n/\ell$, with all these sets pairwise disjoint. For a hyperedge $h \in \cH$, let $W_h = \bigcup_{x\in h}C_x$. We form $\Gamma$ on $\bigcup_{x\in V(\cH)}C_x$ as follows. Independently for each $h\in\cH$ and $x\in h$, choose $i_{h,x}\in[\ell]$ uniformly at random. Conditional on these choices, choose all values of the maps $\phi_h:W_h\to V(G)$ independently, with $\phi_h(v)$ uniform in $V_{i_{h,x}}$ for $v\in C_x$. For distinct $x,y\in V(\cH)$ and vertices $v\in C_x$, $w\in C_y$, put an edge between $v$ and $w$ if and only if there is a hyperedge $h$ containing $x,y$ such that $\phi_h(v)$ and $\phi_h(w)$ are adjacent in $G$. Put no edges inside any $C_x$.

Let $m=v(\Gamma)=nr^2/\ell$. Regardless of the random choices, we have $\Delta(\Gamma)\le drn/\ell=dm/r$. By the property \eqref{prop:transversal}, the sets $\bigcup_{x\in X_i}C_x$ give an $r$-equicolouring of $\Gamma$. By linearity, the sets $W_h$ corresponding to distinct hyperedges intersect in at most one block $C_x$. Since each $C_x$ is independent, edges contributed by other hyperedges do not change the graph induced on $W_h$. Consequently, for every $h\in\cH$ and distinct $v,w\in W_h$, we have $vw\in E(\Gamma)$ if and only if $\phi_h(v)\phi_h(w)\in E(G)$. We show that, with positive probability, for every $S\subseteq V(\Gamma)$ of size at least $\beta m$ there exists $h\in\cH$ with $|\phi_h(S\cap W_h)|\ge(1-\alpha)n$.

Fix $A\subseteq V(\cH)$ of size $|A| \ge \beta r^2/2$, and let $E(A)$ be the set of hyperedges satisfying~\eqref{eq:robust-rich-lines}, so $|E(A)|\ge dr/2$. Call $h\in E(A)$ \emph{balanced} if every $i\in[\ell]$ occurs among the choices $i_{h,x}$ for $x\in A\cap h$ at least $|A|/(4 \ell r)$ times. By Chernoff's inequality and a union bound, the probability that $h$ is not balanced is at most
$$
    \ell \exp \left(-\frac{|A|}{16r\ell}\right).
$$
These events are independent for different hyperedges. Thus the probability that fewer than $dr/4$ hyperedges in $E(A)$ are balanced is at most
$$
    2^{dr}\left(\ell e^{-|A|/(16r\ell)}\right)^{dr/4}
    \le e^{-d|A|/(128\ell)}
    \le e^{-\beta dr^2/(256\ell)},
$$
where the first inequality follows from \eqref{eq:alpha_beta_l} and $C$ being sufficiently large.

Consider now some $S\subseteq V(\Gamma)$ with $|S|\ge\beta m$, and set
$$
    A=\{x\in V(\cH) \; :\; |S \cap C_x| \ge \beta n/(2\ell)\}.
$$
At most $\beta m/2$ vertices of $S$ lie in blocks indexed outside $A$. Since each block has $n/\ell$ vertices, it follows that $|A| n / \ell \ge \beta m / 2$, and so $|A|\ge \beta r^2 / 2$. Condition on any fixed choice of all $i_{h,x}$ for which at least $dr/4$ hyperedges $h \in E(A)$ are balanced. For every balanced $h\in E(A)$ and every $i\in[\ell]$, at least
$$
    \frac{\beta n}{2\ell} \cdot \frac{|A|}{4\ell r} 
$$
vertices of $S\cap W_h$ have their images sampled from $V_i$. Hence, for any $T'\subseteq V(G)$ with $|T'|\ge\alpha n$, conditional independence gives
$$
    \Pr\bigl[\phi_h(S\cap W_h)\cap T'=\emptyset\bigr]
    \le \exp\left(-\frac{\beta |A|}{8r\ell}\sum_{i=1}^{\ell}|T'\cap V_i|\right)
    \le e^{-\alpha\beta n |A|/(8r\ell)}.
$$
Taking a union bound over $T'$ gives
$$
    \Pr\bigl[|\phi_h(S\cap W_h)|<(1-\alpha)n\bigr]
    \le 2^n e^{-\alpha\beta n|A|/(8r\ell)}
    \le e^{-\alpha\beta n|A|/(16r\ell)},
$$
which follows from \eqref{eq:alpha_beta_l} and $C$ being sufficiently large. The maps on distinct hyperedges are conditionally independent, so the probability that all balanced hyperedges fail is at most
$$
    e^{-\alpha\beta nd |A|/(64\ell)} \le e^{-\alpha\beta^2dm/128}.
$$
This bound is uniform over the conditioned choices of $i_{h,x}$. 

Taking a union bound over the at most $2^{r^2}$ sets $A$ bounds the probability that some $A$ has too few balanced hyperedges. Conditional on any choice of the labels for which this does not happen, a union bound over the at most $2^m$ sets $S$ bounds the probability that some $S$ has no successful hyperedge. Thus the total failure probability is at most
$$
  \exp\left(r^2 - \frac{\beta dr^2}{256\ell} \right) + \exp\left(m - \frac{\alpha\beta^2 dm}{128}\right) < e^{-r^2}+e^{-m}<1,
$$
for sufficiently large $C$.

Fix a successful choice of $\Gamma$. For a given $S\subseteq V(\Gamma)$ of size at least $\beta m$, choose $h\in\cH$ with $|\phi_h(S\cap W_h)|\ge(1-\alpha)n$, and set $T=\phi_h(S\cap W_h)$. Choosing one preimage in $S\cap W_h$ of each vertex of $T$ gives an induced copy of $G[T]$ in $\Gamma[S]$. Since $|T|\ge(1-\alpha)n$, the assumption on $G$ implies that $G[T]\ramseyind(F,H)$ for every graph $F$ on at most $b$ vertices. Hence $\Gamma$ is $(H,b,\beta)$-Ramsey. 
\end{proof}

\section{Joining subgraphs via the Sauer--Shelah inequality}
\label{sec:join}

Given graphs $J\subseteq\Gamma$ on the same vertex set, we write $F\subseteq(J,\Gamma)$ if $\Gamma$ contains an induced copy of $F$ all of whose edges lie in $J$. We refer to such a copy as \emph{a copy of $F$ in $(J,\Gamma)$}. For $b\in\mathbb N$, we say that $(J,\Gamma)$ is \emph{$b$-universal} if $F\subseteq(J,\Gamma)$ for every graph $F$ on at most $b$ vertices. The following lemma is the central part of the proof of Theorem~\ref{thm:main}.

\begin{lemma}\label{lemma:join}
  There exists an absolute constant $C>1$ such that the following holds. Let $0 < q \le 1/4$, and let $\Gamma_1$ and $\Gamma_2$ be graphs on disjoint vertex sets $V_1$ and $V_2$, respectively, with $|V_1|=|V_2|=m$ and maximum degree at most $m^{49/50}$, for some integer $m \ge 1$. Let $b \ge 1$ be an integer satisfying
  $$
    m \ge q^{-Cb}.
  $$
  Then there exists a graph $\Gamma$ on $V_1 \cup V_2$ with $\Gamma[V_i]=\Gamma_i$ for $i\in\{1,2\}$, such that the following hold for all sets $X_i\subseteq V_i$ with $|X_i| \ge m^{99/100}$:
  \begin{enumerate}[(i)]
    \item\label{join:density} $e_\Gamma(X_1,X_2)\ge(1-2q)|X_1||X_2|$.

    \item\label{join:copies} Let $J \subseteq \Gamma[X_1 \cup X_2]$ be such that $e_J(X_1,X_2)\ge q |X_1||X_2|$, and $(J[U],\Gamma[U])$ is $b$-universal for each $i \in \{1,2\}$ and every $U\subseteq X_i$ with $|U|\ge m^{49/50}$. Then $(J,\Gamma[X_1\cup X_2])$ is $2b$-universal.
  \end{enumerate}
\end{lemma}

A related extension statement appears in Arag\~ao et al.~\cite[Lemma 3.1]{AragaoCamposDahiaFilipeMarciano2025}; their main proof actually uses the stronger version \cite[Lemma 5.1]{AragaoCamposDahiaFilipeMarciano2025}. Their proof uses hypergraph containers to extend copies by one vertex. Here we join two families of copies using dependent random choice and the Sauer--Shelah inequality. We first recall these two ingredients. The first is a result of Conlon, Fox, and Sudakov \cite[Lemma~2.1]{ConlonFoxSudakov2013}, proved using dependent random choice.

\begin{lemma} \label{lemma:drc}
Let $b,s\in\mathbb{N}$, $h>0$, and $0<\delta\le1$. Let $B$ be a bipartite graph with nonempty parts $X_1$ and $X_2$, and suppose $B$ has at least $\delta |X_1||X_2|$ edges. If
$$
 \binom{|X_2|}{b}\left(\frac{h}{|X_1|}\right)^s \le\frac{\delta^s |X_2|}{2},
$$
then there exists $U_2 \subseteq X_2$ with $|U_2| \ge \delta^s |X_2|/2$, such that every $b$-subset of $U_2$ has at least $h$ common neighbours in $X_1$.
\end{lemma}

The second ingredient is the classical Sauer--Shelah inequality~\cite{sauer72,shelah72}. Given a set $V$ and a family $\cF \subseteq 2^V$, we say that $\cF$ \emph{shatters} a set $S \subseteq V$ if $\{F\cap S:F\in\cF\}=2^S$.

\begin{lemma} \label{lemma:shatter}
    Suppose $V$ is a set of size $n$. If a family $\cF \subseteq 2^V$ does not shatter any set of size $d$, then
    $$
        \left| \cF \right| \le \sum_{i = 0}^{d-1} \binom{n}{i}.
    $$
\end{lemma}

We use the following probabilistic consequence of Lemma~\ref{lemma:shatter}. A related weighted double-counting argument appears in~\cite[Lemma~3.6]{BBDFP2025}.

\begin{lemma} \label{lem:weighted-sauer-shelah}
Let $\Omega$ be a finite set, and let $0 \le d \le |\Omega|$ be an integer. For every $N \subseteq \Omega$, let
$$
    \mathcal{J}_N \subseteq 2^{\Omega\setminus N}
$$
be a family closed under taking subsets. For each $S \subseteq \Omega$, define
$$
    \mathcal{A}_S := \bigl\{ J_X \subseteq S \; : \; J_X \in \mathcal J_{S \setminus J_X} \bigr\}.
$$
Suppose that no family $\mathcal A_S$ shatters a set of size $d+1$. Let $0<q<1$, and let $N \subseteq \Omega$ be a random set obtained by including each element of $\Omega$ independently with probability $q$. Then, for any $k \ge 0$,
$$
    \Pr\bigl[ \exists J_X \in \mathcal{J}_N \; : \; |J_X| \ge k \bigr] \le (1-q)^k \sum_{i=0}^{d}\binom{|\Omega|}{i} q^i.
$$
\end{lemma}
\begin{proof}
By Lemma~\ref{lemma:shatter}, for every $S \subseteq \Omega$ we have
\begin{equation}
\label{eq:shatter_S}
    |\mathcal A_S| \le \sum_{i=0}^{d}\binom{|S|}{i},
\end{equation}
where $\binom{|S|}{i}=0$ when $i > |S|$. Define the weighted count of $\mathcal{J}_N$ as
$$
    Z(N) := \sum_{J_X \in\mathcal J_N} \left(\frac{q}{1-q}\right)^{|J_X|}.
$$
The weights are chosen so that the contribution of each pair $(N,J_X)$ to $\mathbb E Z(N)$ depends only on its union $N\cup J_X$, allowing us to apply~\eqref{eq:shatter_S}. If some $J_X \in \mathcal J_N$ has size at least $k$, then every subset of $J_X$ also belongs to $\mathcal J_N$, so
\begin{equation}
\label{eq:large_J}
    Z(N) \ge \sum_{J' \subseteq J_X} \left(\frac{q}{1-q}\right)^{|J'|} = \left(1+\frac{q}{1-q}\right)^{|J_X|} = (1-q)^{-|J_X|} \ge (1-q)^{-k}.
\end{equation}
We next bound $\mathbb E Z(N)$. Expanding the expectation, we have
\begin{align*}
    \mathbb E Z(N)  &= \sum_{N\subseteq\Omega} q^{|N|}(1-q)^{|\Omega|-|N|} \sum_{J_X\in\mathcal J_N} \left(\frac{q}{1-q}\right)^{|J_X|} \\
    &= \sum_{N\subseteq\Omega} \sum_{J_X\in\mathcal J_N} q^{|N|+|J_X|} (1-q)^{|\Omega|-|N|-|J_X|}.
\end{align*}
Every pair appearing in this sum satisfies $J_X\cap N=\varnothing$.
Writing $S=N\cup J_X$, we have $N=S\setminus J_X$, and the condition
$J_X\in\mathcal J_N$ becomes $J_X\in\mathcal A_S$. Thus, regrouping the terms according to $S$ yields
$$
    \mathbb E Z(N) = \sum_{S\subseteq\Omega} q^{|S|}(1-q)^{|\Omega|-|S|} |\mathcal A_S|.
$$
Applying \eqref{eq:shatter_S}, we obtain
$$
    \mathbb E Z(N) \le \sum_{i=0}^{d} \sum_{S\subseteq\Omega} q^{|S|}(1-q)^{|\Omega|-|S|} \binom{|S|}{i}.
$$
For each $i$, the inner sum is the expected number of $i$-element subsets contained in a random $q$-subset of $\Omega$. Each of the $\binom{|\Omega|}{i}$ possible subsets is contained in that random set with probability $q^i$, thus, by the linearity of expectation, we have
$$
    \sum_{S\subseteq\Omega} q^{|S|}(1-q)^{|\Omega|-|S|} \binom{|S|}{i} = \binom{|\Omega|}{i}q^i.
$$
Finally, \eqref{eq:large_J} and Markov's inequality imply
$$
    \Pr\bigl[ \exists J_X\in\mathcal J_N:\ |J_X|\ge k \bigr] \le (1-q)^k \, \mathbb E Z(N) \le (1-q)^k \sum_{i=0}^{d}\binom{|\Omega|}{i}q^i.
$$
\end{proof}

Let us now briefly discuss the proof strategy. We obtain $\Gamma[X_1, X_2]$ by sampling each crossing edge with probability $1-q$. Fix a target $F$ on $2b$ vertices, sets $X_1,X_2$, and an internal subgraph $J'$ satisfying the $b$-universality assumption. Write $\Omega=X_1\times X_2$, and let $N$ be the random set of absent crossing pairs. Note that each element is included in $N$ with probability $q$, independently. The crossing edge sets $J_X\subseteq\Omega\setminus N$ that fail to complete a copy of $F$ form a family $\mathcal J_N$ closed under taking subsets. We apply Lemma~\ref{lem:weighted-sauer-shelah} to bound the probability that this family contains a set of size at least $q|\Omega|$.

The key is to bound the sizes of sets shattered by the associated families $\mathcal A_S$. Dependent random choice and universality of $J'$ ensure that every $B\subseteq\Omega$ with $|B|>q^4|\Omega|$ contains all crossing pairs between suitable copies of $F[W_1]$ and $F[W_2]$, for any equipartition $V(F)=W_1\cup W_2$. Shattering $B$ by a family $\mathcal{A}_S$ would realise the prescribed pattern of edges and non-edges between these copies by some $J_X \in \mathcal{A}_S$, completing a copy of $F$ and thus contradicting $J_X \in \mathcal{J}_{S \setminus J_X}$. Lemma~\ref{lem:weighted-sauer-shelah} then gives a failure probability at most $\exp\{-q^3|\Omega|\}$. This permits a union bound over $F$, $X_1,X_2$, and $J'$, with the number of choices of $J'$ controlled by the maximum degree assumption.

\begin{proof}[Proof of Lemma \ref{lemma:join}]
We obtain $\Gamma$ by adding each edge between $V_1$ and $V_2$ with probability $1 - q$, independently, and show that both properties hold simultaneously with positive probability. By Chernoff's inequality and a union bound, the probability that the first property fails is easily seen to be less than $1/2$, with room to spare; we omit the routine calculation. We now focus on the second property.

Fix $X_i\subseteq V_i$ with $|X_i|\ge m^{99/100}$, and let $J' \subseteq \Gamma[X_1]\cup\Gamma[X_2]$ be such that $(J'[U],\Gamma[U])$ is $b$-universal for each $i\in\{1,2\}$ and every $U\subseteq X_i$ with $|U|\ge m^{49/50}$. Furthermore, fix a graph $F$ on $2b$ vertices. We show that
\begin{equation} \label{eq:union_bound_pr}
\Pr\left[
    \begin{aligned}
        \exists J_X \subseteq \Gamma[X_1, X_2] \;:\;&
        \begin{gathered}
        e_{J_X}(X_1,X_2)\ge q|X_1||X_2| \text{ and }\\
        F\not\subseteq (J' \cup J_X, \Gamma[X_1 \cup X_2])
        \end{gathered}
    \end{aligned}
\right] \le \exp\{-q^3 |X_1||X_2| \}.
\end{equation}
Let us first explain why this bound suffices. For fixed $X_1$ and $X_2$, the degree assumption on $\Gamma_i$ gives
$$
 e(\Gamma[X_1])+e(\Gamma[X_2]) \le \frac{m^{49/50} (|X_1|+|X_2|)}{2} \le m^{-1/100}|X_1||X_2|.
$$
Therefore, there are at most $2^{m^{-1/100}|X_1||X_2|}$ choices of the edge set of the subgraph $J'$, which determines the event in~\eqref{eq:union_bound_pr}, and at most $2^{4b^2}$ labelled graphs $F$ with $2b$ vertices. For sufficiently large $C$, the assumption $m \ge q^{-Cb}$ ensures
$$
 m^{-1/100}|X_1||X_2| + 4b^2 < \frac{1}{2} q^3 |X_1||X_2|.
$$
Therefore, a union bound using~\eqref{eq:union_bound_pr} shows that property \eqref{join:copies}, restricted to graphs with $2b$ vertices, fails for fixed sets $X_1$ and $X_2$ with probability at most $\exp\{-q^3 |X_1||X_2| / 2\}$. There are at most $4^m$ choices of the sets $X_1$ and $X_2$. By taking $C$ to be sufficiently large, the probability that there exist choices of $X_1$ and $X_2$, each of size at least $m^{99/100}$, for which \eqref{join:copies} fails is therefore at most
$$
 4^m e^{-q^3 m^{99/50} / 2} < 1/2.
$$
Together with the bound for the first property, this gives a positive probability that both properties hold. A successful choice of $\Gamma$ gives every target on $2b$ vertices. For a smaller target, first add isolated vertices and then restrict the resulting copy. Thus~\eqref{eq:union_bound_pr} implies the lemma.

We now turn to the proof of \eqref{eq:union_bound_pr}. Let $V(F) = W_1 \cup W_2$ be an arbitrary partition with $|W_1| = |W_2| = b$, and set $F_i = F[W_i]$. For brevity, let $\Omega = X_1 \times X_2$ denote the set of all possible edges between $X_1$ and $X_2$. 
Given disjoint sets $J_X,N\subseteq\Omega$, we say that $(J_X,N)$ is \emph{avoiding} if there are no copies of $F_i$ in $(J'[X_i],\Gamma[X_i])$, for $i \in \{1,2\}$, whose embeddings together map every edge of $F$ between $W_1$ and $W_2$ into $J_X$ and every non-edge between these parts into $N$. For a set $N \subseteq \Omega$, set
$$
    \mathcal{J}_N = \{ J_X \subseteq \Omega \setminus N \; : \; (J_X, N) \text{ is avoiding}\}.
$$
Let $N=\Omega\setminus E(\Gamma[X_1,X_2])$, so that each pair belongs to $N$ independently with probability $q$. The event on the left-hand side of~\eqref{eq:union_bound_pr} implies that $\mathcal{J}_N$ contains a set of size at least $q|\Omega|$. It therefore suffices to show that
$$
 \Pr\bigl[\exists J_X\in\mathcal{J}_N:\ |J_X|\ge q|\Omega|\bigr]
 \le \exp\{-q^3|\Omega|\}.
$$

The family $\mathcal{J}_N$ is closed under taking subsets. To apply Lemma~\ref{lem:weighted-sauer-shelah}, we bound the sizes of sets shattered by the associated families
$$
 \mathcal A_S=\{J_X\subseteq S:(J_X,S\setminus J_X)\text{ is avoiding}\}.
$$ 

\begin{claim}\label{claim:joining-copies}
Let $B \subseteq \Omega$ satisfy $|B| > q^4 |\Omega|$. Then there exist copies of $F_1$ and $F_2$ in $(J'[X_1], \Gamma[X_1])$ and $(J'[X_2], \Gamma[X_2])$, respectively, on vertex sets $\Phi_i \subseteq X_i$ such that $\Phi_1 \times \Phi_2 \subseteq B$.
\end{claim}
\claimproofstart
With Lemma \ref{lemma:drc} in mind, set $s = 200b$, $h = m^{49/50}$, and $\delta = q^4$. For sufficiently large $C$, a routine calculation shows
$$
    \binom{|X_2|}{b} \left( \frac{h}{|X_1|} \right)^{s} < h < \frac{\delta^s |X_2|}{2}.
$$
Thus Lemma \ref{lemma:drc} gives a set $X_2' \subseteq X_2$ with $|X'_2| \ge h$ such that every $b$-subset of $X_2'$ has at least $h$ common neighbours in $X_1$. By the assumption on $J'$, we have $F_2 \subseteq (J'[X_2'], \Gamma[X_2'])$; denote the vertex set of such a copy by $\Phi_2$. Since $|\Phi_2| = b$, its common neighbourhood $X_1' \subseteq X_1$ in $B$ has size at least $h$. Applying the assumption on $J'$ gives a copy of $F_1$ in $(J'[X_1'], \Gamma[X_1'])$ on a vertex set $\Phi_1$. Then $\Phi_1 \times \Phi_2 \subseteq B$, as desired.
\claimproofend

The previous claim implies that, for any set $S \subseteq \Omega$, the family $\mathcal{A}_S$ does not shatter a set larger than $q^4 |\Omega|$. This can be seen as follows. Suppose, towards a contradiction, that $\mathcal{A}_S$ shatters a set $B\subseteq S$ with $|B|>q^4|\Omega|$. Let $\Phi_1$ and $\Phi_2$ be vertex sets of copies of $F_1$ and $F_2$, as given by Claim \ref{claim:joining-copies}, so that $\Phi_1\times\Phi_2\subseteq B$. Since $\mathcal{A}_S$ shatters $B$, we can choose $J_X\in\mathcal{A}_S$ whose intersection with $\Phi_1\times\Phi_2$ consists precisely of the pairs corresponding to edges of $F$ under these embeddings. The pairs corresponding to non-edges of $F$ then lie in $S\setminus J_X$, contradicting that $(J_X,S\setminus J_X)$ is avoiding. 

Apply Lemma~\ref{lem:weighted-sauer-shelah} with $d=\lfloor q^4|\Omega|\rfloor$ and $k=q|\Omega|$. By the binomial theorem,
$$
 \sum_{i=0}^{d}\binom{|\Omega|}{i}q^i = \sum_{i = 0}^d \binom{|\Omega|}{i} \frac{q^{4i}}{q^{3i}} \le q^{-3d} \sum_{i = 0}^{|\Omega|} \binom{|\Omega|}{i} q^{4i} = q^{-3d} (1 + q^4)^{|\Omega|}. 
$$
Consequently,
$$
 \Pr\bigl[\exists J_X\in\mathcal J_N:\ |J_X|\ge q|\Omega|\bigr] \le (1-q)^{q|\Omega|}\sum_{i=0}^{d}\binom{|\Omega|}{i}q^i 
 \le \left( (1 - q)^q q^{-3 q^4} (1 + q^4) \right)^{|\Omega|}.
$$
Using $1 + x \le e^x$ for $x \in \mathbb{R}$, and $\ln(1/x) \le 1/x - 1$ for $x > 0$, we further upper bound the preceding expression as
$$
    \exp\{ \left( -q^2 + 3q^4 \ln(1/q) + q^4 \right) |\Omega| \} \le e^{-q^3 |\Omega|},
$$
which follows from $q \le 1/4$. This proves~\eqref{eq:union_bound_pr} and completes the proof.
\end{proof}

\section{Proof of Theorem \ref{thm:main}}

The following is a standard result on embedding graphs in hosts with sufficiently uniform edge distribution. We include a short proof for the sake of completeness.

\begin{lemma}\label{lemma:density_embedding}
Let $H$ be a graph on the vertex set $[t]$, let $p\in(0,1]$, and let $z\in\mathbb{N}$. Suppose that $\Gamma$ is a graph whose vertex set is partitioned into sets $W_1,\ldots,W_t$, each of size at least $tzp^{-t}$. Suppose further that, for every $\{i,j\}\in E(H)$ and all subsets $X_i\subseteq W_i$ and $X_j\subseteq W_j$ with $|X_i|,|X_j|\ge z$, we have
$$
  e_\Gamma(X_i,X_j)\ge |X_i||X_j| p.
$$
Then $\Gamma$ contains a copy of $H$ in which each vertex $i\in[t]$ is mapped to a vertex of $W_i$.
\end{lemma}
\begin{proof}
Embed the vertices in the order $1,\ldots,t$, maintaining a candidate set in each unused part. When embedding vertex $i$, require its image to have at least a $p$-fraction of each future candidate set corresponding to a neighbour of $i$. For each such set, fewer than $z$ vertices fail this requirement: otherwise these vertices and the candidate set would contradict the density assumption. After at most $t-1$ restrictions, every candidate set still has size at least $p^{t-1}tzp^{-t}\ge tz$. The exclusions remove fewer than $tz$ vertices, so there is a choice for the image of $i$. Restrict the candidate set of each future neighbour of $i$ to the neighbours of the chosen image. This completes the induction.
\end{proof}

For a graph $H$ on $t\ge1$ vertices, we abbreviate $(H,b,1-1/(2t))$-Ramsey to $(H,b)$-Ramsey. We now prove the inductive step of Theorem~\ref{thm:main}. Starting with an $(H,b)$-Ramsey graph $G$, we use Lemma~\ref{lemma:robust} to obtain an $(H,b,\beta)$-Ramsey graph $\Gamma$ of controlled order and maximum degree, for a suitably small $\beta>0$. We then take disjoint copies of $\Gamma$, indexed by $V(H)$, and join copies corresponding to edges of $H$ using Lemma~\ref{lemma:join}. In a colouring with no red induced copy of $H$, Lemma~\ref{lemma:density_embedding} forces a pair of sufficiently large vertex sets $X_1,X_2$ with many blue edges. The amplified robustness allows us to apply property~\eqref{join:copies} to these sets, giving a blue induced copy of every graph on at most $2b$ vertices. Deleting at most a $1/(2t)$-fraction of the vertices of the new host leaves at least half the vertices in each copy of $\Gamma$, so the same argument also establishes the required robustness.

The second inequality in the following lemma controls the number of colour classes relative to the order of the new graph. This is needed because the next application of Lemma 2.2 depends on $n'$ as well as on $\ell'$.

\begin{lemma}\label{lemma:doubling}
There exists an absolute constant $C>1$ such that the following holds. Let $H$ be a graph on $t\ge2$ vertices, and let $1\le b\le t$ be an integer. Suppose that $G$ is an $\ell$-equicolourable $(H,b)$-Ramsey graph on $n$ vertices. Then there exist an integer $\ell'\ge t$ and an $\ell'$-equicolourable $(H,2b)$-Ramsey graph $G'$ on $n'$ vertices such that
$$
    n' \le (n\ell)^{6/5} \left(\frac{128t}{b}\right)^{Cb} \qquad\text{and}\qquad \frac{n'}{(\ell')^2}\ge\frac{n}{t\ell}.
$$
\end{lemma}
\begin{proof}
We start by setting up parameters and proving that they satisfy certain inequalities. Let $C_0,C_1$ be the constants in Lemmas~\ref{lemma:robust} and~\ref{lemma:join}, respectively, and choose a sufficiently large absolute constant $D$ in terms of $C_0,C_1$. Set $q=b/(128t)$, and let $r$ be the least power of two such that $r \ge \ell(n\ell)^{1/10}q^{-Db}$. For 
$$
    m = \frac{r^2}{\ell} \, n,
$$
we then have
\begin{equation}\label{eq:host-order}
 (n\ell)^{6/5}q^{-2Db} \le m < 4(n\ell)^{6/5}q^{-2Db}.
\end{equation}
A routine but rather tedious calculation shows that, for $\alpha=1/(2t)$ and $\beta=m^{-1/50}$, we have
\begin{equation} \label{eq:r_d}
    d := \left\lceil C_0 \frac{\ell \log(2 \ell)}{\alpha \beta^2} \right\rceil \le m^{-1/50} r. 
\end{equation}
This is needed to apply Lemma \ref{lemma:robust} with these parameters.

\paragraph{Graph construction.} We apply Lemma~\ref{lemma:robust} to $G$ with $\alpha, \beta$, and the prime power $r$ chosen above. Inequality \eqref{eq:r_d} shows that this is indeed a valid choice for $r$. This gives us an $r$-equicolourable $(H,b,m^{-1/50})$-Ramsey graph $\Gamma$ with
$$
    v(\Gamma) = m \quad \text{and} \quad \Delta(\Gamma) \le \frac{d}{r}m \stackrel{\eqref{eq:r_d}}{\le} m^{49/50}.
$$
Therefore, every set $U\subseteq V(\Gamma)$ with $|U|\ge m^{49/50}$ satisfies
\begin{equation}\label{eq:small-ramsey}
 \Gamma[U]\ramseyind(F,H) \quad \text{for every graph $F$ on at most $b$ vertices.}
\end{equation}

Take disjoint copies $\Gamma_1,\ldots,\Gamma_t$ of $\Gamma$, indexed by the vertices of $H$. By~\eqref{eq:host-order}, we have $m\ge q^{-C_1b}$. For each $ij\in E(H)$, apply Lemma~\ref{lemma:join} to $\Gamma_i,\Gamma_j$ with parameters $m,q,b$, and use the resulting edges between these two parts. Put no edges between $\Gamma_i,\Gamma_j$ when $ij\notin E(H)$. This gives a graph $G'$ on $n'=tm=tnr^2/\ell$ vertices. By ~\eqref{eq:host-order} and $q^{-b}\ge128t$, where the latter is easily derived from $1 \le b \le t$ and $q = b / (128 t)$, we have
$$
    n'=tm<4t(n\ell)^{6/5}q^{-2Db} \le (n\ell)^{6/5}q^{-Cb}
$$
for $C = 2D + 1$. Using a distinct set of $r$ colours on each $\Gamma_i$ gives an $\ell'$-equicolouring of $G'$, where $\ell'=tr\ge t$, and hence
$$
    \frac{n'}{(\ell')^2}=\frac{tnr^2/\ell}{t^2r^2}=\frac{n}{t\ell}.
$$
This proves that $n'$ and $\ell'$ satisfy the desired properties.

\paragraph{Ramsey property.} Consider a set $S\subseteq V(G')$ of size $|S| \ge (1-1/(2t))n'$, and colour the edges of $G'[S]$ red and blue without a red induced copy of $H$. Since at most $m/2$ vertices lie outside $S$, each set $W_i=S\cap V(\Gamma_i)$ has size at least $m/2$. By~\eqref{eq:small-ramsey}, every subset of $W_i$ of size at least $m^{49/50}$ contains a blue induced copy of every graph on at most $b$ vertices.

We claim that there are $ij\in E(H)$ and sets $X_i\subseteq W_i$, $X_j\subseteq W_j$ of size $z = \lceil m^{99/100} \rceil$ with at least $qz^2$ blue edges between them. Otherwise, condition~\eqref{join:density} of Lemma~\ref{lemma:join} gives red density at least $p:=1-3q$ between every such pair of sets. By averaging, the same lower bound holds between any subsets of $W_i,W_j$ of size at least $z$. Moreover,
$$
 |W_i| \ge m/2\ge tzp^{-t},
$$
which follows from \eqref{eq:host-order}, $4tp^{-t} \le 4te^{b/32 }\le m^{1/100}$ (which holds for sufficiently large $D$), and $z \le 2m^{99/100}$. Lemma~\ref{lemma:density_embedding}, applied to the red graph, therefore gives a copy of $H$ meeting each $W_i$ once. It is induced in $G'$ because no edges were placed between parts corresponding to nonedges of $H$, contradicting our choice of colouring.

Let $J$ be the subgraph of $G'[X_i\cup X_j]$ consisting of all its blue edges. By~\eqref{eq:small-ramsey}, the pair $(J[U],G'[U])$ is $b$-universal for every $U\subseteq X_i$ or $U\subseteq X_j$ with $|U|\ge m^{49/50}$. Applying condition~\eqref{join:copies} of Lemma~\ref{lemma:join}, we obtain a blue induced copy of every graph on at most $2b$ vertices. Thus $G'$ is $(H,2b)$-Ramsey.
\end{proof}

\begin{proof}[Proof of Theorem~\ref{thm:main}]
The case $s=1$ is immediate, so assume $2\le s\le t$. Let $C$ be the constant in Lemma~\ref{lemma:doubling}. We show by induction that, for every $1\le b \le 2t$ which is a power of $2$, there exists an $\ell_b$-equicolourable $(H,b)$-Ramsey graph $G_b$ on $n_b$ vertices, where
$$
    n_b \le \left(\frac{C't}{b}\right)^{5Cb}
    \quad \text{ and } \quad
    \ell_b \le \min \left\{\sqrt{n_b}, \frac{n_b}{t}\right\},
$$
for an absolute constant $C'>1$.  This suffices to prove the theorem: take $b$ to be the smallest power of $2$ at least $s$. Since $s\le b<2s$, the definition of $(H,b)$-Ramsey gives
$$
  \log\rind(F,H)\le\log n_b \le 10Cs\log\frac{C't}{s}.
$$

For $b=1$, take $G_1$ to be the empty graph on $t$ vertices. It is trivially $1$-equicolourable and $(H,1)$-Ramsey, and satisfies both bounds. Consider now some $1\le b\le t$ which is a power of $2$, and suppose the claim holds for $b$. We show that it then also holds for $2b$. Apply Lemma~\ref{lemma:doubling} to $G_b$ to obtain an $\ell'$-equicolourable $(H,2b)$-Ramsey graph $G'$ on $n'$ vertices. Then
$$
    \frac{n'}{(\ell')^2 }\ge \frac{n_b}{t\ell_b} \ge1
    \qquad \text{and} \qquad
    \frac{n'}{\ell'} \ge \ell'\ge t,
$$
which gives the required bound on $\ell'$. Moreover, using $\ell_b \le \sqrt{n_b}$, the lemma gives
\begin{equation} \label{eq:inductive_dependence}
    n'\le n_b^{9/5}\left(\frac{128t}{b}\right)^{Cb}.
\end{equation}
The inductive bound on $n_b$ therefore yields
$$
n' \le\left(\frac{C't}{b}\right)^{9Cb} \left(\frac{128t}{b}\right)^{Cb} = \left(\frac{128}{C'}\right)^{Cb} \left(\frac{C't}{b}\right)^{10Cb} \le \left(\frac{C't}{2b}\right)^{5C(2b)},
$$
where the last inequality holds provided $C'$ is sufficiently large. Taking $G_{2b} = G'$ completes the induction. The subquadratic dependence on $n_b$ in~\eqref{eq:inductive_dependence} is what allows the induction to close.
\end{proof}

\bibliographystyle{abbrv}
\bibliography{references}

\end{document}